# Necessary and sufficient conditions for the inverse problem of one class ordinary differential operators with complex periodic coefficients

Efendiev R.F.
Institute of Applied Mathematics BSU,
Z.Khalilova 23, 370148, Baku, Azerbaijan.

rakibaz@yahoo.com

Abstract :
The basic purpose of the present paper is the full solutions of the inverse problem (i.e. a finding of necessary and sufficient conditions) for the operator with complex periodic coefficients.

1. <u>Statement of the problem and the formulation of the basic results</u>.

Let's consider a class $Q^2$ of all $2\pi$ periodic complex valid functions on the real axis $\mathbb{R}$, belonging to space $L_2[0,2\pi]$, and its subclass $Q_+^2$ consisting of the functions of a type

$$p_\gamma(x) = \sum_{n=1}^{\infty} p_{\gamma n} \exp(inx) \ . \ \sum_{\gamma=0}^{2m-2} \sum_{n=1}^{\infty} n^\gamma |p_{\gamma n}| < \infty. \tag{1.1}$$

The basic purpose of the present paper is the full solutions of the inverse problem (i.e. a finding of necessary and sufficient conditions) for the operator L, generated by differential expression

$$l(y) = (-1)^m y^{(2m)} + \sum_{\gamma=0}^{2m-2} p_\gamma(x) y^{(\gamma)}(x) \tag{1.2}$$

in the space $L_2(-\infty,\infty)$, with potentials of a type (1.1).

The inverse problem, for potentials of a type (1.1) for the first time is put and solved in work [1] where is shown, that the equation $l(y) = \lambda^{2m} y$ has the solutions

$$\varphi(x,\lambda\omega_\tau) = e^{i\lambda\omega_\tau x} + \sum_{j=1}^{2m-1} \sum_{\alpha=1}^{\infty} \sum_{n=1}^{\alpha} \frac{V_{n\alpha}^{(j)}}{n + \lambda\omega_\tau(1-\omega_j)} e^{(i\lambda\omega_\tau + i\alpha)x}, \tau = 0, 2m-1, \omega_j = \exp(ij\pi/m) \tag{1.3}$$

and Wronckian of the systems of solutions $\varphi(x,\lambda\omega_\tau)$ being equal to $(i\lambda)^{m(2m-1)} A$, where

$$A = \begin{vmatrix} 1 & 1 & .. & 1 \\ \omega_1 & \omega_2 & .. & \omega_{2m-1} \\ .. & .. & .. & .. \\ \omega_1^{2m-1} & \omega_2^{2m-1} & .. & \omega_{2m-1}^{2m-1} \end{vmatrix},$$

it is not zero at $\lambda \neq 0$, and a limit

$$\varphi_{nj}(x) \equiv \lim_{\lambda \to -\lambda_{nj}} (\lambda + \lambda_{nj})\varphi(x,\lambda). \ n \in N, j = \overline{1, 2m-1}$$

also is the solution of the equation $l(y) = \lambda^{2m} y$, but linearly dependent on $\varphi(x, \lambda_{nj}\omega_j)$, because there are such numbers $\tilde{S}_{nj}, n \in N, j = \overline{1, 2m-1}$, that satisfy the relation

$$\varphi_{nj}(x) = \tilde{S}_{nj}\varphi(x, \lambda_{nj}\omega_j). \ \lambda_{nj} = -\frac{n}{1-\omega_j}. \tag{1.4}$$

Further, M.G.Gasymov [1] has established, that if

I. $\sum\limits_{n=1}^{\infty} n|\tilde{S}_n| < \infty$

II. $4^{m-1} a_m \sum\limits_{n=1}^{\infty} \frac{|\tilde{S}_n|}{n+1} = p < 1$

where

$$a_m = \max_{\substack{1 \le j \le l \le 2m-1 \\ 1 \le n, r < \infty}} \frac{|(1-\omega_j)(n+r)|}{|r(1-\omega_j) - n(1-\omega_l)\omega_j|} \ , \ \tilde{S}_n = \sum_{j=1}^{2m-1} n^{2m-2} |\tilde{S}_{nj}|, \tag{1.5}$$

then exist unique determined functions $p_\gamma(x), \gamma = \overline{0, 2m-2}$ of the form (1.1), for which numbers $\{\tilde{S}_n\}$ are under formulas (1.3) - (1.4).

The full solution of this problem at m=1 is given in work [5](the present paper doesn't include the concerning reference, the detail reference reader can find in [5,7]) where it is proved the following statement:

In order to the set of sequence of complex numbers $\{S_n\}$ was a set of the spectral data of the operator $L = \left(\dfrac{d}{dx}\right)^2 + p_0(x)$, with potential $p_0(x) \in Q_+^2$, it is necessary and sufficient the simultaneously satisfying of the conditions:

1) $\{nS_n\}_{n=1}^{\infty} \in l_1$;

2) an infinite determinant

$$D(z) \equiv \left\| \delta_{nk} + \frac{2S_k}{n+k} e^{i\frac{n+k}{2}z} \right\|_{n,k=1}^{\infty}$$

exists (here and further we use the definition $\delta_{rn}$ is Kronecker delta, $E_n$ identity matrix of order $n \times n$), it is continuous, does not zero in closed half plane $\overline{C_+} = \{z : \text{Im } z \ge 0\}$ and analytical inside of open half plane $C_+ = \{z : \text{Im } z > 0\}$.

In the present work, using a technique of works [2], [5] the full inverse problem is solved for the operator (1.2) with potentials of the type (1.1). Now let's formulate the basic result of the present work.

<u>Definition</u>. The sequence $\{\tilde{S}_{nj}\}_{n=1, j=1}^{\infty, 2m-1}$ constructed with the help of formulas (1.4), we call as a set of the spectral data of the operator L generated by differential expression (1.2), with potential (1.1).

<u>Thoerem1</u>. In order to the given sequence of complex numbers $\{\widetilde{S}_{nj}\}_{n=1,j=1}^{\infty,2m-1}$ would be the set of the spectral data of operator L generated by differential expression (1.2) and potential (1.1), it is necessary and sufficient that the conditions satisfied in the same time;

1) $\{n\widetilde{S}_n\}_{n=1}^{\infty} \in l_1$; (1.6)

2) The infinite determinant

$$D(z) \equiv \det\left\|\delta_{rn}E_{2m-1} + \left\|\frac{i(1-\omega_l)\widetilde{S}_{nj}}{r\omega_l(1-\omega_j)-n(1-\omega_l)}e^{i\frac{n}{1-\omega_j}z}e^{-i\frac{r\omega_l}{1-\omega_l}z}\right\|_{j,l=1}^{2m-1}\right\|_{r,n=1}^{\infty} \quad (1.7)$$

exists, is continuous, is not zero in closed half plane $\overline{C_+}=\{z:\operatorname{Im} z \geq 0\}$ and analytical inside of open half plane $C_+ = \{z:\operatorname{Im} z > 0\}$.

<u>2. About the inverse problem of the theory of scattering on half-line</u>.

Taking $x = it, \lambda = -ik, y(x) = Y(t)$ in the equation $l(y) = \lambda^{2m}y$, we receive

$$(-1)^m Y^{(2m)}(t) + \sum_{\gamma=0}^{2m-2} Q_\gamma(t) Y^{(\gamma)}(t) = k^{2m} Y(t), \quad (2.1)$$

where

$$Q_\gamma(t) = (-1)^m (-i)^\gamma \sum_{n=1}^{\infty} p_{\gamma n} e^{-nt}, \quad \sum_{\gamma=0}^{2m-2}\sum_{n=1}^{\infty} n^\gamma |p_{\gamma n}| < \infty. \quad (2.2)$$

As a result we have the equation (2.1) where potential exponentially decreases $t \to \infty$.

<u>3. About the operator of transformation.</u>

In [1] it is proved, that (2.1) has solutions type of

$$f(x, k\omega_\tau) = e^{ik\omega_\tau x} + \sum_{j=1}^{2m-1}\sum_{\alpha=1}^{\infty}\sum_{n=1}^{\alpha} \frac{V_{n\alpha}^{(j)}}{in + k\omega_\tau(1-\omega_j)}e^{(ik\omega_\tau - \alpha)x}, \quad \tau = \overline{0, 2m-1} \quad (3.1)$$

where numbers $V_{n\alpha}^{(j)}$ are determined from the following recurrent formulas:

$$\left[\left(\alpha - \frac{n}{(1-\omega_j)}\right)^{2m} - \left(\frac{n}{(1-\omega_j)}\right)^{2m}\right]V_{n\alpha}^{(j)} = (-1)^{m+1}\sum_{\gamma=0}^{2m-1}\sum_{s=n}^{\alpha-1}\left[i\left(s-\frac{n}{(1-\omega_j)}\right)\right]^\gamma P_{\gamma,s-n}V_{ns}^{(j)} \quad (3.2)$$

at $\alpha = 2, 3, \ldots; n = 1, 2, \ldots, \alpha - 1; j = 1, 2, \ldots, 2m-1$,

$$i^\gamma p_{\gamma\alpha} + \sum_{j=1}^{2m-1}\sum_{n=1}^{\alpha} d_{j\gamma}(n,\alpha)V_{n\alpha}^{(j)} + \sum_{\nu=\gamma+1}^{2m-2}\sum_{j=1}^{2m-1}\sum_{r+s=\alpha}\sum_{n=1}^{s} d_{j\gamma}(n,s,\nu)p_{\nu r}V_{ns}^{(j)} = 0, \quad (3.3)$$

where

$$\frac{1}{n+k(1-\omega_j)}\left[(i\alpha+k)^{2m} - k^{2m} - (i\alpha+k_{nj})^{2m} + k_{nj}^{2m}\right] = \sum_{\gamma=0}^{2m-2} d_{j\gamma}(n,\alpha)k^\gamma; \quad j = \overline{1, 2m-1},$$

$$\frac{(is+k)^\nu - (is+k_{nj})^\nu}{in+k(1-\omega_j)} = \sum_{\gamma=0}^{\nu-1} d_{j\gamma}(n,s,\nu)k^\gamma,$$

and a series (3.1) are 2m time term by term differentiable.

Consequence: Let's the condition (2.2) is satisfied. Then

$$f(t,k) = e^{ikt} + \int_t^\infty K(t,u)e^{iku} du , \qquad (3.4)$$

where

$$K(t,u) = \sum_{j=1}^{2m-1} \sum_{n=1}^{\infty} \sum_{\alpha=n}^{\infty} \frac{V_{n\alpha}^{(j)}}{i(1-\omega_j)} e^{\left[-\alpha t + \frac{n}{1-\omega_j}(t-u)\right]}. \qquad (3.5)$$

## 4. The basic equations of the inverse problem.

Using a technique of work the [1,2] it is possible to receive equality such as (1.4) for solutions (3.1)

$$f_{nj}(t) = S_{nj} f(t, k_{nj}\omega_j), \qquad (4.1)$$

where

$$f_{nj}(t) = \lim_{k \to k_{nj}} [in + k(1-\omega_j)] f(t,k), \quad k_{nj} = -\frac{in}{1-\omega_j}.$$

Equality (4.1) we shall write as

$$\sum_{\alpha=n}^{\infty} V_{n\alpha}^{(j)} e^{-\alpha t} e^{\frac{n}{1-\omega_j}t} = S_{nj} e^{\frac{n\omega_j}{1-\omega_j}t} + \sum_{l=1}^{2m-1} \sum_{r=1}^{\infty} \sum_{\alpha=r}^{\infty} \frac{i(1-\omega_j) V_{nr}^{(l)} S_{nj}}{n\omega_j(1-\omega_l) - r(1-\omega_j)} e^{(-\alpha + \frac{n\omega_j}{1-\omega_j})t} \qquad (4.2)$$

Let's multiply both parts of (4.2) by $\frac{1}{i(1-\omega_j)} e^{-\frac{n}{1-\omega_j}u}$ and denote

$$\tilde{F}(t+u) = \sum_{j=1}^{2m-1} \sum_{n=1}^{\infty} \frac{S_{nj}}{i(1-\omega_j)} e^{\frac{n}{1-\omega_j}(t\omega_j - u)}, t \le u , \qquad (4.3)$$

then from the equation (4.2) we shall receive Marchenko's type equations

$$K(t,u) = \tilde{F}(t+u) + \int_t^\infty K(t,s)\tilde{F}(s+u)ds . \qquad (4.4)$$

Lemma 1: If the coefficients $Q_\gamma(t)$ of the equation (2.1) are (2.2) type, then at all $t \ge 0$ the operator of transformation (3.5) satisfies Marchenko's type equation (4.4) in which the function of transition $\tilde{F}(t)$ has type (4.3), and numbers $S_{nj}$ are determined by equality (4.1). From this it is obtained that $S_{nj} = V_{nn}^{(j)}$.

## 4.1 Resolvability of the basic equation and uniqueness of the solution of the inverse problem.

By nucleus of the operator of transformation the coefficients $Q_\gamma(t)$ are reconstructed with the help of recurrent formulas (3.2) - (3.3). As consequences, the basic equation (4.4) and type (4.3) of the functions of transition of the inverse problem put natural statement about reconstructed of coefficients of the equation (2.1) on numbers $S_{nj}$. In this statement is an important moment the proof of unequivocal resolvability of the basic equation (4.4).

<u>Lemma 2</u>. The homogeneous equation

$$g(s) - \int_0^\infty \widetilde{F}(u+s)g(u)du = 0 \tag{4.5}$$

corresponding to potentials $Q_\gamma(t) \in Q_+^2$ has only trivial solution.

The proof: Let $g \in L_2(R^+)$ be solution of the equation (4.5) and let $f$ be solution of

$$f(s) + \int_0^s K(t,s)f(t)dt = g(s) \tag{4.6}$$

Substituting $g$ in (4.6) and taking into account the equation (4.4) we receive

$$f(s) + \int_0^s K(t,s)f(t)dt + \int_0^\infty [f(u) + \int_0^u K(t,u)f(t)dt]\widetilde{F}(u+s)du =$$

$$= f(s) + \int_s^\infty f(t)[\widetilde{F}(t+s) + \int_t^\infty K(t,u)\widetilde{F}(u+s)du] = 0.$$

As at all $t \geq s$ the estimation

$$\left|\widetilde{F}(t+s) + \int_0^\infty K(t,u)\widetilde{F}(u+s)du\right| \leq Ce^{-s}$$

is satisfied, it follows, that $f = 0, g = 0$ and lemma it is proved. From this lemma follows

<u>Theorem 2:</u> Coefficients $Q_\gamma(t)$ of the equation (2.1) satisfying condition (2.2) is unequivocally determined by numbers $S_{nj}$.

In the equation (1.2) we shall replace $x$ by $x+a$, where $\operatorname{Im} a \geq 0$. Then we shall receive equation of the same type with potential $Q_\gamma^a(x) = Q_\gamma(x+a)$ satisfying the condition (1.1). Note, that the functions $\varphi(x+a, \lambda\omega_j)$ are solutions of the equation

$$(-1)^m y^{(2m)}(x) + \sum_{\gamma=0}^{2m-2} Q_\gamma^a(x) y^{(\gamma)}(x) = \lambda^{2m} y(x)$$

which, at $x \to \infty$ has a form

$$\varphi(x+a, \lambda\omega_j) = e^{i\lambda\omega_j a} e^{i\lambda\omega_j x} + o(1).$$

Therefore functions

$$\varphi^a(x, \lambda\omega_j) = e^{-i\lambda\omega_j a} \varphi(x+a, \lambda\omega_j)$$

are solutions of type (1.3). We shall denote further through $S_{nj}(a)$ the spectral data of operator L with potential $Q_\gamma^a(x)$

$$L \equiv (-1)^m \frac{d^{2m}}{dx^{2m}} + \sum_{\gamma=0}^{2m-2} Q_\gamma^a(x) \frac{d^\gamma}{dx^\gamma}.$$

According to (1.4) we have

$$S_{nj}(a) \varphi^a(x, \lambda_{nj}\omega_j) = \lim_{\lambda \to \lambda_{nj}}[n + \lambda(1-\omega_j)]\varphi^a(x,\lambda) = \lim_{\lambda \to \lambda_{nj}}[n + \lambda(1-\omega_j)]\varphi(x+a,\lambda)e^{-i\lambda a}$$

$$= e^{i\frac{n}{1-\omega_j}a} S_{nj} \varphi(x+a, \lambda_{nj}\omega_j) e^{i\frac{n}{1-\omega_j}a} e^{-i\frac{n\omega_j}{1-\omega_j}a} = S_{nj} \varphi^a(x, \lambda_{nj}\omega_j) e^{ina} = S_{nj} \varphi^a(x, \lambda_{nj}\omega_j),$$

as consequence

$$S_{nj}(a) = e^{ina} S_{nj} \tag{4.7}$$

Now considering as above we receive basic equations of type (4.4) with the transitive function

$$\tilde{F}_a(t+u) = \sum_{j=1}^{2m-1} \sum_{n=1}^{\infty} \frac{S_{nj}(a)}{i(1-\omega_j)} e^{\frac{n}{1-\omega_j}(t\omega_j - u)} = \tilde{F}(t - ia + u - ia) = \tilde{F}(t + u - 2ia)$$

and validity of the following lemma.

<u>Lemma 3</u>: At each fixed value $a, (\operatorname{Im} a \geq 0)$ the homogeneous equation

$$g(s) - \int_0^\infty \tilde{F}(u+s-2ia)g(u)du = 0 \tag{4.8}$$

has only trivial solution in the space $L_2(R^+)$.

<u>5. The proof of the basic theorem 2.</u>

<u>Necessity:</u> From the relation (4.1) and a type of function $f_{nj}(t)$ we shall receive $S_{nj} = V_{nn}^j$.

Therefore

$$\sum_{j=1}^{2m-1} \sum_{n=1}^{\infty} n^{2m-1} |S_{nj}| \leq \sum_{j=1}^{2m-1} \sum_{n=1}^{\infty} n^{2m-1} |V_{nn}^j| < \infty,$$

i.e. $n^{2m-1} |S_{nj}| \in l_1$. Necessity of the condition (1) of the theorem is proved. For the proof of necessity of a condition (2) at first all we shall show, that from trivial resolvability of the basic equation (4.4) at $t = 0$ in a class of functions satisfying to the inequality $\|g(u)\| \leq Ce^{-\frac{u}{2}}, u \geq 0$, follows trivial resolvability of the infinite systems of the equations in $l_2(1, \infty, R^{2m-1})$

$$g_{jn} - \sum_{l=1}^{2m-1} \sum_{r=1}^{\infty} \frac{i(1-\omega_j) S_{rl}}{n\omega_j(1-\omega_l) - r(1-\omega_j)} g_{lr} = 0 \quad . \tag{5.1}$$

Really, if $\{g_{jn}\} \in l_2, j = \overline{1, 2m-1},$ and the solution of this system exists, then the function

$$g(u) = (g_1(u), g_2(u), \ldots g_{2m-1}(u)) = \sum_{j=1}^{2m-1} \sum_{n=1}^{\infty} S_{nj} g_{jn} e^{-\frac{n}{1-\omega_j} u} \tag{5.2}$$

is determined for all $u \geq 0$ and satisfies inequality

$$|g(u)| \leq Ce^{-\frac{u}{1-\omega_j}}; u \geq 0.$$

also is the solution of the equation (4.5), because of

$$g(u) - \int_0^{\infty} g(s) \widetilde{F}(u+s) ds = \sum_{j=1}^{2m-1} \sum_{n=1}^{\infty} S_{nj} g_{jn} e^{-\frac{n}{1-\omega_j} u} - \int_0^{\infty} (\sum_{j=1}^{2m-1} \sum_{n=1}^{\infty} S_{nj} g_{jn} e^{-\frac{n}{1-\omega_j} s})(\sum_{l=1}^{2m-1} \sum_{r=1}^{\infty} \frac{S_{rl}}{i(1-\omega_j)} e^{-\frac{r}{1-\omega_l}(s\omega_l - u)}) ds =$$

$$= \sum_{j=1}^{2m-1} \sum_{n=1}^{\infty} S_{nj} g_{jn} e^{-\frac{n}{1-\omega_j} u} - \sum_{j=1}^{2m-1} \sum_{n=1}^{\infty} \sum_{l=1}^{2m-1} \sum_{r=1}^{\infty} \frac{S_{nj} S_{rl}}{i(1-\omega_l)} g_{jn} e^{-\frac{r}{1-\omega_l} u} \int_0^{\infty} e^{-\frac{n}{1-\omega_j} s} e^{-\frac{r\omega_l}{1-\omega_l} s} ds =$$

$$= \sum_{j=1}^{2m-1} \sum_{n=1}^{\infty} S_{nj} g_{jn} e^{-\frac{n}{1-\omega_j} u} - \sum_{j=1}^{2m-1} \sum_{n=1}^{\infty} \sum_{l=1}^{2m-1} \sum_{r=1}^{\infty} \frac{i(1-\omega_j) S_{nj} S_{rl}}{n\omega_j(1-\omega_l) - r(1-\omega_j)} g_{lr} e^{-\frac{n}{1-\omega_j} u} =$$

$$= \sum_{j=1}^{2m-1} \sum_{n=1}^{\infty} S_{nj} e^{-\frac{n}{1-\omega_j} u} [g_{jn} - \sum_{l=1}^{2m-1} \sum_{r=1}^{\infty} \frac{i(1-\omega_j) S_{rl}}{n\omega_j(1-\omega_l) - r(1-\omega_j)} g_{lr}] = 0.$$

Therefore, $g(u) = 0$, $S_{nj} g_{jn} = 0$ at all, $n \geq 1$, $j = \overline{1,2m-1}$, hence $g_{jn} = 0$, $j = \overline{1,2m-1}$, $n \geq 0$ by (5.1). Now we shall introduce in space $l_2(1,\infty, R^{2m-1})$ the operator $F(t)$ defined by matrix

$$F_{rn}(t) = \left\|F_{rn}^{jl}\right\|_{j,l=1}^{2m-1} = \left\|\frac{i(1-\omega_l)S_{nj}}{r\omega_l(1-\omega_j)-n(1-\omega_l)} e^{-\frac{n}{1-\omega_j}t} e^{\frac{r\omega_l}{1-\omega_l}t}\right\|_{j,l=1}^{2m-1} \quad (5.3)$$

and let $\varphi_{2k-1} = \{(\delta_{v1})\delta_{rk}\}_{v,r=1}^{2m-1,\infty}$, $\varphi_{2k} = \{(\delta_{v2})\delta_{rk}\}_{v,r=1}^{2m-1,\infty}$ $\{(\delta_{v1})$-column vector$\}$ ortonormal system in this space, as $n^{2m-1}|S_{nj}| \in l_1$, we have

$\sum_{j,k=1}^{\infty} \left|(F\varphi_j, \varphi_k)_{l_2(1,\infty, R^{2m-1})}\right| < \infty$, and the operator $F(t)$ is nuclear [3]. Therefore there exists a determinant $\Delta(t) = \det(E - F(t))$ of the operator $E - F(t)$ connected as it is easy to see, with a determinant from the condition 2) of the Theorem2 by the relations $\Delta(-iz) = \det(E - F(-iz)) \equiv D(z)$. Determinant of system (5.1) is $D(0)$, and a determinant of the similar system considering to potential $Q_\gamma^z = Q_\gamma(x+z)$, $\text{Im}\, z \geq 0$ is

$$D(z) \equiv \det\left\|\delta_{rn} E_{2m-1} - \left\|\frac{i(1-\omega_l)S_{nj}(z)}{r\omega_l(1-\omega_j)-n(1-\omega_l)}\right\|_{j,l=1}^{2m-1}\right\|_{r,n=1}^{\infty} =$$

$$= \det\left\|\delta_{rn} E_{2m-1} - \left\|\frac{i(1-\omega_l)S_{nj}}{r\omega_l(1-\omega_j)-n(1-\omega_l)} e^{inz}\right\|_{j,l=1}^{2m-1}\right\|_{r,n=1}^{\infty} =$$

$$= \det\left\|\delta_{rn} E_{2m-1} - \left\|\frac{i(1-\omega_l)S_{nj}}{r\omega_l(1-\omega_j)-n(1-\omega_l)} e^{i\frac{n}{1-\omega_j}z} e^{-i\frac{r\omega_l}{1-\omega_l}z}\right\|_{j,l=1}^{2m-1}\right\|_{r,n=1}^{\infty}.$$

Therefore for the proof of necessity of the condition 2) of the theorem should be checked up, that $\Delta(0) = D(0) \neq 0$. The system (5.1) can be written as the equation in $l_2(1,\infty, R^{2m-1})$,

$$g - F(0)g = 0.$$

Since the operator $F(0)$ is nuclear, we can apply to this equation Fredholm theory, according to which, its trivial resolvability is equivalent to that $\det(E - F(0))$ is not zero [6]. Necessity of a condition 2) is proved.

<u>Sufficiency:</u> Let's multiply the equation (4.4) $e^{\frac{r\omega_l}{1-\omega_l}u}$ and integrate on $u \in [t,\infty)$. Then we receive

$$k(t) = F(t)e(t) + k(t)F(t) \quad (5.4)$$

in which the operator $F(t)$ is defined by the matrix $\|F_{rn}(t)\|_{r,n=1}^{\infty}$ of a type (5.3),

$$e(t) = \|e_{nj}(t) = e^{\frac{n\omega_j}{1-\omega_j}t}\|_{n,j=1}^{\infty,2m-1} \quad ; \quad k(t) = \|k_{lr}(t)\|_{l,r=1}^{2m-1,\infty} = \left\|\int_t^{\infty} K(t,u) e^{\frac{r\omega_l}{1-\omega_l}u} du\right\|_{l,r=1}^{2m-1,\infty} ;$$

As $F(t)$ is nuclear at $t \geq 0$ and conditions $\Delta(t) = \det(E - F(t)) \neq 0$ holds, there exists bounded in $l_2$ the inverse operator $R(t) = (E - F(t))^{-1}$. Considering $F(t)e(t) \in l_2$, from (5.4) we find that

$$k(t) = R(t)F(t)e(t) \tag{5.5}$$

Defining $\langle f, g \rangle = \sum_{n=1}^{\infty} f_n g_n$, from (4.4) receive

$$K(t,u) = \langle e(t), A(u) \rangle + \langle k(t), A(u) \rangle = \langle e(t), A(u) \rangle + \langle R(t)F(t)e(t), A(u) \rangle = \langle e(t) + R(t)F(t)e(t), A(u) \rangle =$$

$$= \langle R(t)e(t), A(u) \rangle, \tag{5.6}$$

where $A(u)$ is defined by the matrix

$$A(u) = \|a_{jn}(u)\|_{j,n=1}^{2m-1,\infty} = \left\|\frac{S_{nj}}{i(1-\omega_j)} e^{-\frac{n}{1-\omega_j}u}\right\|_{j,n=1}^{2m-1,\infty}$$

Let's assume, now, that conditions of the theorem are satisfied. According to the stated reasons we define at $0 \leq t \leq u$ the function $K(t,u)$ by equality (5.6). Then at $u \geq t$ we have

$$K(t,u) - \int_t^{\infty} K(t,s)\widetilde{F}(s+u)ds = \langle R(t)e(t), A(u) \rangle - \int_t^{\infty} \langle R(t)e(t), A(s)\langle e(s), A(u) \rangle \rangle ds =$$

$$= \langle R(t)e(t), A(u) \rangle - \left\langle R(t)e(t), \left\langle \int_t^{\infty} A(s)e(s)ds, A(u) \right\rangle \right\rangle = \langle R(t)e(t), A(u) \rangle -$$

$$- \langle R(t)e(t), \langle A(u), F(t) \rangle \rangle = \langle R(t)e(t), A(u) \rangle - \langle R(t)e(t), A(u)F^*(t) \rangle = \langle R(t)e(t), A(u) - A(u)F^*(t) \rangle =$$

$$= \langle e(t), A(u) \rangle = \widetilde{F}(t+u) \tag{5.7}$$

where "*" designates transition to the matrix connected with $F(t)$ rather bilinear form $\langle .,. \rangle$. The following lemma is proved.

<u>Lemma 4:</u> At everyone $t \geq 0$ the nucleus $K(t,u)$ of the operator of transformation satisfies to the basic equation

$$K(t,u) = F(t+u) + \int_t^{\infty} K(t,s)\widetilde{F}(s+u)ds$$

where $\widetilde{F}(t+u)$ it is defined from (4.3).
Resolvability of the basic equation follows from the lemma 2 unequivocal. By substitution it is easy to calculate, that the solution of the basic equation is

$$K(t,u) = \sum_{j=1}^{2m-1} \sum_{n=1}^{\infty} \sum_{\alpha=n}^{\infty} \frac{V_{n\alpha}^{(j)}}{i(1-\omega_j)} e^{\left[-\alpha t + \frac{n}{1-\omega_j}(t-u)\right]},$$

where numbers $V_{n\alpha}^{(j)}$ are determined from the recurrent relations

$$V_{nn}^{(j)} = S_{nj},$$

$$V_{n\alpha+n}^{(j)} = (1-\omega_j)V_{nn}^{(j)} \sum_{l=1}^{2m-1}\sum_{r=1}^{\alpha} \frac{V_{r\alpha}^{(l)}}{r(1-\omega_j)-n(1-\omega_l)\omega_j}.$$

Passing to the proof of the basic statement that coefficients $Q_\gamma(t)$ look like (2.2), all over again we shall establish for the matrix elements $R_{rn}(t)$ of the operator $R(t)$ the estimation

$$\left|R_{rn}^{jl}(t)\right| \leq \delta_{rn}\delta_{jl} + CS_n; \tag{5.8}$$

$$\left|\frac{\partial^{2m-\tau}}{\partial t^{2m-\tau}} R_{rn}^{jl}(t)\right| \leq CS_n; \tau = \overline{1,2m-1} \tag{5.9}$$

where $C = \max\{C_k > 0, k = \overline{1,2m-1}\}$ - a constant, and $S_n = \sum_{j=1}^{2m-1} n^{2m-1}\left|S_{nj}\right|$;

Really from equality $R(t) = E + R(t)F(t)$ follows

$$\left|R_{rn}^{lj}(t)\right| \leq \delta_{rn}\delta_{lj} + \sum_{\tau=1}^{2m-1}\sum_{p=1}^{\infty}\left|R_{rp}^{l\tau}(t)\right|\left|F_{pn}^{\tau j}(t)\right| \leq \delta_{rn}\delta_{lj} + 2\sum_{\tau=1}^{2m-1}(\sum_{p=1}^{\infty}\left|R_{rp}^{l\tau}(t)\right|^2)^{\frac{1}{2}}(\sum_{p=1}^{\infty}\left|F_{pn}^{\tau j}(t)\right|^2)^{\frac{1}{2}} \leq \delta_{rn}\delta_{lj} +$$

$$+ C_1 a_m ((R(t)R^*(t))_{pp} \sum_{p=1}^{\infty}\frac{1}{(n+p)^2})S_n \leq \delta_{rn}\delta_{lj} + C_1\|R(t)\|_{l_2 \to l_2} S_n.$$

On the other hand, as it has already been shown, operator - function $R(t) = (E - F(t))^{-1}$ exists and is bounded in $l_2$ (valid nuclear $F(t)$ at $t \geq 0$ and conditions $\Delta(t) = \det(E - F(t)) \neq 0$), that proves the first inequality (5.8). The proof of the second estimation (5.9), we do similarly, proceeding from identity $\frac{d^l}{dt^l}R(t) = \sum_{n=1}^{l} C_l^n (\frac{d^{l-n}}{dt^{l-n}}R(t))(\frac{d^n}{dt^n}F(t))R(t), l = \overline{1,2m-1}$, and using the first estimation (5.8). Then

$$\left|\frac{d}{dt}R_{rn}^{lj}(t)\right| \leq \sum_{\tau=1}^{2m-1}\sum_{p,q=1}^{\infty}\left|R_{rp}^{l\tau}(t)\right|\left|\frac{d}{dt}F_{pq}^{\tau k}(t)\right|\left|R_{qn}^{kj}(t)\right| \leq \sum_{\tau=1}^{2m-1}\sum_{p,q=1}^{\infty}(\delta_{rp}\delta_{l\tau}+C_2 S_p)S_q(\delta_{qn}\delta_{kj}+C_3 S_n) \leq$$

$$\leq (1 + C_4 \sum_{p=1}^{\infty} S_p)^2 S_n \leq CS_n.$$

Further, proceeding from equality

$$\frac{d^l}{dt^l}R(t) = \sum_{n=1}^{l} C_l^n (\frac{d^{l-n}}{dt^{l-n}}R(t))(\frac{d^n}{dt^n}F(t))R(t), l = \overline{1,2m-1},$$

with the help of a mathematical induction method the inequality (5.9) is proved. In [4] it is established validity of the following relations (for conformity with our case, we shall put $q_{2m-2-\gamma}(x) = Q_\gamma(x)$).

$$(-1)^m \frac{\partial^{2m}}{\partial x^{2m}}K(x,t) + \sum_{\gamma=0}^{2m-2} q_{2m-2-\gamma}(x)\frac{\partial^\gamma}{\partial x^\gamma}K(x,t) - \frac{\partial^{2m}}{\partial t^{2m}}K(x,t) = 0$$

$$q_0(x) = 2m\frac{d}{dx}K(x,x)$$

$$q_{k+1}(x) = \sum_{\nu=0}^{k} q_\nu(x)\sum_{s=\nu}^{k} C_{n-3-s}^{k-s}\left\{\frac{\partial^{s-\nu}}{\partial x^{s-\nu}}K(x,t)\bigg|_{t=x}\right\}^{(k-s)} +$$

$$+ \sum_{k=0}^{k+2} C_{n-1-\nu}^{k+2-\nu} \left\{ \frac{\partial^\nu}{\partial x^\nu} K(x,t) \bigg|_{t=x} \right\}^{(x+2-\nu)} - (-1)^k \frac{\partial^{k+2}}{\partial t^{k+2}} K(x,t) \bigg|_{t=x}$$

$k = 0,1,\ldots 2m - 3$.

Further, it is easy to show, that

$$q_0(x) = \sum_{j=1}^{2m-1} \sum_{n=1}^{\infty} \frac{n \cdot S_{nj}}{i(1-\omega_j)} e^{-nt} + \Pi_0(t)$$

where

$$\Pi_0(t) = \sum_{n,p,q,\tau=1}^{\infty} R_{np}^{ej} F_{pq}^{j\tau} e_{q\tau} e_{\tau q} = <R(t)F(t)e(t), A(t)>$$

is $2i\pi$ periodic function and has the bounded derivative up to the order (2m-1). Then Fourier coefficients $\sum_{n=1}^{\infty} \left| n^{2m-1} \Pi_n \right|^2 < \infty$. But then $\sum_{n=1}^{\infty} n^{2m-2} |\Pi_n| < \infty$. Thus Fourier coefficients $P_{2m-2,n}$ of the function $Q_{2m-2}(x) = q_0(x)$ satisfy to the condition (2.2). Similarly, for all other coefficients $Q_\gamma(x)$, $\gamma = \overline{0,2m-3}$ it is established, that the Fourier coefficients $p_{\gamma n}$ functions $Q_\gamma(x) = q_{2n-2-\gamma}(x)$, $\gamma = \overline{0,2m-3}$ satisfy a condition (2.2). Thus Fourier coefficients of function $P_\gamma(x)$, $\gamma = \overline{0,2m-2}$ satisfied to the condition (1.1).

Let, at least $\{\widetilde{S}_{nj}\}$ be the set of spectral data of the operator $(L - k^{2m}E)$ with constructed coefficients $P_\gamma(x)$. For the end the proof it needs to be shown, that $\{S_{nj}\}$ coincides with an initial set $\{\widetilde{S}_{nj}\}$. It turns out from equality $\widetilde{S}_{nj} = V_{nn}^{(j)} = S_{nj}$.
The theorem is proved.

### Acknowledgments
The author has benefited from discussions with M.G.Casymov and I.M.Guseynov.